\documentclass[12pt]{amsart}
\usepackage{amsmath,amsthm,amsfonts,amssymb,amscd}
\usepackage{fancyhdr}
\usepackage{setspace}
\newtheorem{problem}{PROBLEM}
\newtheorem{theorem}{THEOREM}[section]
\theoremstyle{definition}
\newtheorem{corollary}{Corollary}[theorem]
\newtheorem{definition}[theorem]{Definition}
\newtheorem*{remark}{Remark}
\newtheorem{example}{Example}[section]
\newtheorem{lemma}[theorem]{Lemma}
\newtheorem{proposition}[theorem]{Proposition}
\numberwithin{equation}{section}
\newcommand{\C}{\mathbb{C}}

\newcommand{\x}{\zeta}

\newcommand{\by}{\overline{\y}}
\newcommand{\bx}{\overline{\x}}
\renewcommand{\a}{\alpha}
\renewcommand{\b}{\beta}
\renewcommand{\r}{\gamma}
\newcommand{\s}{\sigma}
\newcommand{\de}{\delta}

\newcommand{\y}{\eta}

\newcommand{\bt}{\begin{theorem}}
\newcommand{\et}{\end{theorem}}
\newcommand{\bco}{\begin{corollary}}
\newcommand{\eco}{\end{corollary}}
\newcommand{\bd}{\begin{definition}}
\newcommand{\ed}{\end{definition}}
\newcommand{\bp}{\begin{problem}}
\newcommand{\ep}{\end{problem}}
\newcommand{\bl}{\begin{lemma}}
\newcommand{\el}{\end{lemma}}
\newcommand{\bprop}{\begin{proposition}}
\newcommand{\eprop}{\end{proposition}}
\newcommand{\br}{\begin{remark}}
\newcommand{\er}{\end{remark}}
\newcommand{\bpf}{\begin{proof}}
\newcommand{\epf}{\end{proof}}
\title{Iterated Antiderivative Extensions}
\author{V. Ravi Srinivasan}
\date{\today}
\address{\noindent Department of Mathematics and Computer Science, Rutgers University, Newark, NJ 07102.}
\email{ravisri@rutgers.edu}
\begin{document}
\maketitle

\begin{abstract} Let $F$ be a characteristic zero differential field with an algebraically closed field of constants and  let $E$ be a no new constants extension of $F$. We say that $E$ is an \textsl{iterated antiderivative extension} of $F$ if $E$ is a liouvillian extension of $F$ obtained by adjoining antiderivatives alone. In this article, we will show that if $E$ is an iterated antiderivative extension of $F$ and $K$ is a differential subfield of $E$  that contains $F$ then $K$ is an iterated antiderivative extension of $F$.

\end{abstract}
\section{Introduction}

Let $F:=\C(z)$ be the differential field of rational functions in one complex variable $z$ with the usual derivation $d/dz$. Consider the liouvillian extensions $E_1:=F(e^{z^2},\int e^{z^2})$ and $E_2:=F(\sqrt{1-z^2},$ $\sin^{-1}z)$ of $F$. In \cite{Rosenlicht-Singer}, M. Rosenlicht and M. Singer shows that the differential subfield $F((\int e^{z^2})/e^{z^2})$ of $E_1$ and the differential subfield $F(\sqrt{1-z^2}\sin^{-1} z)$ of $E_2$ are not liouvillian extensions of $F$. Thus, differential subfields of liouvillian extensions, in general, need not be liouvillian. However, if $ L:=\C(z,\log z,\log(\log z))$ then one can list all the differential subfields of $ L$ that contains $\C$ and they are $\C$, $\C(z)$, $\C(z,\log z)$ and $ L$, see example \ref{log example}. Clearly, in this case, all the differential subfields are liouvillian. Thus, it is of considerable interest to know when differential subfields of a liouvillian extension are liouvillian? In this article, we will show that if a liouvillian extension is obtained by adjoining antiderivatives alone then  its differential subfields can also be obtained by adjoining antiderivatives alone. This is the main result of this article and it appears  as theorem \ref{structure IAE}. An analogue of theorem \ref{structure IAE} for generalized elementary extensions can be found in \cite{Rosenlicht-Singer} and \cite{M.Sing}.

\subsection{Differential Fields:} Let $F$ be a field of characteristic zero. A \textsl{derivation} on a field $F$, denoted by $'$, is an additive map $':F\to F$ that satisfies the Leibniz law $(xy)'=x'y+xy'$ for every $x,y\in F$. A field equipped
with a derivation map is called a \textsl{differential field}.  The set of \textsl{constants} of a differential field is the kernel of the map $'$ and it can be seen that the set of constants is a differential subfield of $F$.  Let $E$ and $F$  be differential fields. We say that $E$ is a \textsl{differential field extension}  of $F$ if $E$ is a field extension of $F$ and the restriction of the derivation of $E$ to $F$ coincides with the derivation of $F$. A differential field extension $E$ of $F$ is called a \textsl{no new constants} extension if the constants of $E$ are the same as the constants of $F$.

Throughout this article, we fix a ground differential field $F$ of characteristic zero. All the differential fields considered henceforth are either differential subfields of $F$ or a differential field extension of $F$. We reserve the notation $'$ to denote the derivation map of any given differential field.

 Let $E$ be a no new constants extension of
$F$. An element $\x\in E$ is called an \textsl{antiderivative}
(of an element) of $F$
if $\x'\in F$. We say that
$E$ is an \textsl{antiderivative extension} of $F$ if $E=F(\x_1,\x_2,\cdots,\x_n)$, where $\x_1,\x_2,\cdots,\x_n$ are antiderivatives
of  $F$. Elements $\x_1,\x_2,\cdots,\x_n\in E$ are called \textsl{iterated antiderivatives} of $F$ if $\x'_1\in F$ and for $i\geq 2$, $\x'_i\in F(\x_1,\x_2\cdots,\x_{i-1})$. We call
$E$ an \textsl{iterated antiderivative extension} of $F$ if $E=F(\x_1,\x_2,\cdots,\x_n)$, where  $\x_1,\x_2,\cdots,\x_n$ are iterated antiderivatives of $F$. And, if $E=F(\x_1,\x_2,\cdots,\x_n)$ and for each $i$,  $\x'_i\in F(\x_1,\cdots,\x_{i-1})$ or $\x'_i/\x_i\in F(\x_1,\cdots,\x_{i-1})$ or $\x_i$ is algebraic over $F(\x_1,\cdots,\x_{i-1})$ then we call $E$ a \textsl{liouvillian} extension of $F$. Now it is clear that the differential fields $E_1, E_2$ and $ L$, mentioned in the beginning of this article, are examples of liouvillian extensions of $\C$ and that $ L$ is an iterated antiderivative extension of $\C$. A field automorphism of $E$ that fixes the elements of $F$ and commutes with the derivation is called a differential field automorphism and the group of all such automorphisms will be denoted by $G(E|F)$. That is, $G(E|F)=\{\s\in Aut(E|F)|\s(y)'=\s(y')\ \ \text{for all} \ \ y\in E\}$.

 Every antiderivative extension of  $F$ is an iterated antiderivative extension of $F$. But the converse is not true: for example,  consider the differential field $\mathbb{C}( z,\log  z)$ with the usual derivation $d/dz$, where $\mathbb{C}$ is the field of complex numbers. Clearly, $\mathbb{C}( z,\log  z)$ is an iterated antiderivative extension of $\C$. Observe that all the antiderivatives of the field $\mathbb{C}$ are of the form $c z+ d$ where $c,d\in \mathbb{C}$ and since $\log  z\notin \C( z)$, we see that $\mathbb{C}( z,\log  z)$ is not an antiderivative extension of $\mathbb{C}$.

%%%%%%%%%%%%%%%%%%%%%%%%%%%%%%%%%%%%%%%%%%%%%%%%%%%%%%%%%%%%%%%%%%%%%%%%%%%%%%%%%%%%%%%%%%
%%%%%%%%%%%%%%%%%%%%%%%%%%%%%%%%%%%%%%%%%%%%%%%%%%%%%%%%%%%%%%%%%%%%%%%%%%%%%%%%%%%%%%%%%%
\section{preliminary results}
%%%%%%%%%%%%%%%%%%%%%%%%%%%%%%%%%%%%%%%%%%%%%%%%%%%%%%%%%%%%%%%%%%%%%%%%%%%%%%%%%%%%%%%%%%%%%%%%%%%%%%%%%%%%%
%%%%%%%%%%%%%%%%%%%%%%%%%%%%%%%%%%%%%%%%%%%%%%%%%%%%%%%%%%%%%%%%%%%%%%%%%%%%%%%%%%%%%%%%%%%%%%%%%%%%%%%%%%%%%
It is a well known fact that if $E$ is a no new constants extension of $F$ and if $\x\in E$ is an antiderivative of an element of $F$ then either
$\x$ is transcendental over $F$ or $\x\in F$. Please see \cite{Magid 1} page 7, or \cite{Rosenlicht-Singer} page 329 for a proof. Using this fact,
we will now show that every iterated antiderivative extension of $F$ is a purely transcendental extension of $F$.

%%%%%%%%%%%%%%%%%%%%%%%%%%%%%%%%%%%%%%%%%%%%%%%%%%%%%%%%%%%%%%%%%%%%%%%%%%%%%%%%%%%%%%%%%%%%%%%%%%%%%%%%%%%%%
%%%%%%%%%%%%%%%%%%%%%%%%%%%%%%%%%%%%%%%%%%%%%%%%%%%%%%%%%%%%%%%%%%%%%%%%%%%%%%%%%%%%%%%%%%%%%%%%%%%%%%%%%%%%%%
\bt\label{no alg extn thm}

Let $E$ and $K$ be differential subfields of some no new constants extension of $F$. Suppose that $E=F(\x_1,\x_2,\cdots,\x_n)$ is an iterated antiderivative extension of $F$ and that $K\supseteq F$. Then $KE:=K(\x_1,\x_2,\cdots,\x_n)$ is an iterated antiderivative extension of $K$. Furthermore, If $KE\neq K$ then the set $\{\x_1,\x_2,\cdots,\x_n\}$ contains algebraically independent iterated antiderivatives $\y_1,\y_2,\cdots,\y_t$ of $K$ such that $KE=K(\y_1,\y_2,\cdots,\y_t)$.

\et

\bpf
 Since $K$ contains $F$, it is easy to see that $\x'_i\in K(\x_1,\x_2,\cdots,\x_{i-1})$ and thus
$KE$ is an iterated antiderivative extension of $K$.  Assume that $K(E)\neq K$. To find a transcendence base for $KE$, consisting of iterated antiderivatives of $K$, we use an induction on $n$.
 Case n=1: Since $KE=K(\x_1)\neq K$, we have $\x_1\notin K$. And since $\x'_1\in F\subseteq K$, as noted earlier, $\x_1$ is transcendental over $K$. Set $\y_1:=\x_1$ to prove the theorem. Assume the theorem for n-1 iterated antiderivatives. Induction step: Choose $l$ smallest such that $\x_l\notin K$ and set $\y_1:=\x_l$. Since $\x_1,\cdots,\x_{l-1}\in K$, we see that $\y_1$ is an antiderivative of $K$ and since $\y_1\notin K$, $\y_1$ is transcendental over $K$. Note that $KE$ is generated as a field by $n-l$ iterated antiderivatives, namely $\x_{l+1},\cdots,\x_n$, and the differential field $K(\y_1)$. Now we may apply induction to the iterated antiderivative extension $KE$ of $K(\y_1)$ and obtain iterated antiderivatives $\y_2,\cdots,\y_t\in \{\x_{l+1},\cdots,\x_n\}$ of $K(\y_1)$ such that $\y_2,\cdots,\y_t$ are algebraically independent over $K(\y_1)$ and that $KE=K(\y_1)(\y_2,\cdots,\y_t)$.  \epf
 In theorem \ref{no alg extn thm}, if we choose $K=F$, we obtain that $E$ is a purely transcendental extension of $F$ with a transcendence base consisting of iterated antiderivatives of $F$. Note that theorem \ref{no alg extn thm} is valid for antiderivative extensions as well. Thus, hereafter, when we say $E=F(\x_1,\x_2,\cdots,\x_t)$ is an antiderivative extension or an iterated antiderivative extension of $F$, it is understood that
$\x_1,\x_2,\cdots,\x_t$ are algebraically independent
over $F$. We will use the notation tr.d.$(E|F)$ to denote the transcendence degree of any field extension $E$ over $F$.

\begin{corollary}\label{trd cor}

Let $E$ be an iterated antiderivative extension of $F$ and let $K_1$ and $K$ be differential subfields of $E$. If $K_1\supset K\supseteq F$ then tr.d.$(K_1|F)>$ tr.d.$(K|F)$.

\end{corollary}
\bpf
  Suppose that $K_1\supset K$. Then we have $E\supset K$ and therefore from theorem \ref{no alg extn thm}, we know that $KE=E$ is a purely transcendental  extension of $K$. Thus if $u\in K_1-K$ then $u\in E-K$ and therefore $u$ is transcendental over $K$. Thus tr.d.$(K_1|K)\geq 1$. Note that tr.d.$(K_1|F)$= tr.d.$(K_1|K)+$ tr.d.$(K|F)$ and that tr.d. $K_1|F<\infty$ since tr.d.$(E|F)<\infty$. Hence tr.d.$(K_1|F)>$ tr.d.$(K|F)$.
\epf

Let $M$ be a differential field extension of $F$. We call $M$, a \textsl{minimal} differential field extension of $F$ if $M\supset F$ and if $K$ is a differential subfield of $M$ such that $M\supseteq K\supseteq F$ then $K=M$ or $K=F$.

\begin{corollary}\label{ext min}
Let $E, K$ and $K_1$ be as in corollary \ref{trd cor}. Then $K_1$ contains a minimal differential field extension of $K$.
\end{corollary}

\bpf
If $K_1$ is not a minimal differential field extension of $K$ then it contains a proper subfield $K_1\supset M\supset K$. And,
from corollary \ref{trd cor}, we know that tr.d.$(K_1|K)>$ tr.d.$(M|K)$. Since tr.d.$(K_1|K)<\infty$, the rest of the proof follows by an induction on tr.d.$(K_1|F)$.\epf

%%%%%%%%%%%%%%%%%%%%%%%%%%%%%%%%%%%%%%%%%%%%%%%%%%%%%%%%%%%%%%%%%%%%%%%%%%%%%%%%%%%%%%%%%%%%%%%%%%%%%%%%%%%%%%%%%%%%%%%%%%%%%%%%%%%%%%%%%%%%%%%%%%%%%%%
%%%%%%%%%%%%%%%%%%%%%%%%%%%%%%%%%%%%%%%%%%%%%%%%%%%%%%%%%%%%%%%%%%%%%%%%%%%%%%%%%%%%%%%%%%%%%%%%%%%%%%%%%%%%%%%%%%%%%%%%%%%%%%%%%%%%%%%%%%%%%%%%%%%%%%%

\bt\label{no dlog}

Let $E$  be an iterated antiderivative extension of $F$ and $K\supseteq F$ be a differential subfield of $E$. If there is an element $u\in E$ such  that
$u'/u\in K$ then $u\in K$.

\et

\bpf
To avoid triviality, we may assume $E\neq K$. We observe from theorem \ref{no alg extn thm} that $E=K(\y_1,\y_2,\cdots,\y_t)$ is an iterated antiderivative extension of $K$.
 Let $u\in E$ and $u'/u\in K$. We will use an induction on $t$ to prove our proposition.
 Assume that if $u\in K(\y_1,\y_2,\cdots,\y_{t-1})$ then $u\in K$. Write $u=P/Q$, where $P, Q\in K(\y_1,\y_2,\cdots,\y_{t-1})[\y_t]$ are relatively prime polynomials and $Q$ is monic. Then $u'=(P'Q-Q'P)/Q^2$ and since $f:=u'/u\in K$, we obtain
 $$QPf=P'Q-Q'P.$$ Since $P$ and $Q$ are relatively prime, we then obtain $P$ divides $ P'$ and $Q$
divides $Q'$. Now the facts that, $Q$ is monic, deg
$Q'\lneqq$ deg $Q$ and $Q$ divides $Q'$, all together, will force $Q=1$.
Thus $u=P$ and $P'=fP$. Write $P=\sum^n_{i=0}a_i\y^i_t$ with $a_n\neq 0$ and observe
that $$a'_n\y^n_t+(a'_{n-1}+na_n\y'_t)\y^{n-1}_t+\cdots+a_1\y'_t+a'_0=f(\sum^n_{i=0}a_i\y^i_t)$$ and comparing the leading coefficients, we obtain $a'_n=fa_n .$ Thus  $(u/a_n)'=0$. Since $E$ is a no new constants extension of $F$, there is a
$c\in C$ such that $u=ca_n$.  Now $a_n\in
K(\y_1,\y_2,\cdots,\y_{t-1})$ will imply $u\in
K(\y_1,\y_2,\cdots,\y_{t-1})$.\epf

\begin{remark} Consider the differential field $K:=\C(z,\log z)$ with the derivation $d/dz$, $\overline{K}$ be its algebraic closure and let $u\in \overline{K}-K$. We claim that for any iterated antiderivative extension $E$ of $\C$, the element $u\notin E$. First we note that if $E\neq \C$ is an iterated antiderivative extension of $\C$ with the derivation $d/dz$ then $z\in E$. Now,  suppose that the claim is false. Then by applying \ref{no alg extn thm} to the iterated antiderivative extension $E(\log z)$ of $\C$ we obtain a contradiction to the assumption that $u\notin K$. Thus, if $u=\sqrt{z}+\sqrt[5]{\log z}$, then there are no polynomials $P$, $Q\in \C[x_1,x_2,x_3]$ such that $\sqrt{z}+\sqrt[5]{\log z}=$ $\frac{P(z, \log z, \log(\log z))}{Q(z,\log (z+1), Li_2(z))}$, where $Li_2(z)$ is the dilogarithm $-\int^z_0\frac{\log(1-w)}{w}dw$.

Similarly, as an application of theorem \ref{no dlog}, one can obtain that $e^{\a z}$, where $\a\in \C-\{0\}$ and $e^{-z^2}$ are not in any iterated antiderivative extension of $\C$. In particular, $\int e^{-z^2}$ is not in any iterated antiderivative extension  of $\C$, and thus cannot be expressed in terms of logarithms or polylogarithms.

\end{remark}

%%%%%%%%%%%%%%%%%%%%%%%%%%%%%%%%%%%%%%%%%%%%%%%%%%%%%%%%%%%%%%%%%%%%%%%%%%
%%%%%%%%%%%%%%%%%%%%%%%%%%%%%%%%%%%%%%%%%%%%%%%%%%%%%%%%%%%%%%%%%%%%%%%%%%
\section{structure of antiderivative extensions}
%%%%%%%%%%%%%%%%%%%%%%%%%%%%%

The following theorem characterizes the algebraic dependence of antiderivatives and will be used in numerous occasions in this article.  In this section we will use this theorem to describe the structure of differential subfields of antiderivative extensions.

\bt\label{ostrowski} Let $E\supset F$ be a  no new constants
extension and for $i=1,2,\cdots,n,$ let $\x_i\in E$ be
antiderivatives of $F$. Then either $\x_i$'s are algebraically
independent over $F$ or there is a tuple
$(\a_1,\cdots,\a_n)\in C^n-\{\vec{0}\}$ such that
$\sum^n_{i=1}\a_i\x_i\in F$.

\et

\bpf see \cite{j.ax}, page 260 or \cite{Ravi.S}, page 9. \epf

%%%%%%%%%%%%%%%%%%%%%%%%%%%%%%%%%%%%%%%%%%%%%%%%%%%%%%%%%%%%%%%%%%%%%%%%%%%%%%%%%%%%%%%%%%%%%%%%%%%%%%%%%%%%%
%%%%%%%%%%%%%%%%%%%%%%%%%%%%%%%%%%%%%%%%%%%%%%%%%%%%%%%%%%%%%%%%%%%%%%%%%%%%%%%%%%%%%%%%%%%%%%%%%%%%%%%%%%%%%

\begin{proposition}\label{alg dep of int}
Let $E=F(\x_1,\x_2,\cdots,\x_t)$ be an antiderivative extension of $F$. An element $\x\in E$ is an antiderivative of $F$
if and only if there is a tuple $(\a_1,\cdots,\a_t)\in C^t-\{\vec{0}\}$ and an element $a_\x\in F$ such that $\x=\sum^t_{i=1}\a_i\x_i+a_\x.$
\end{proposition}

\bpf
Let $\x\in E$ be an antiderivative of $F$.
The set $\{\x, \x_1,\x_2,\cdots,\x_t\}$ contains $t+1$ antiderivatives of $F$ and therefore has to be algebraically dependent over $F$. We apply theorem \ref{ostrowski} and obtain constants $\b_i,\r\in C$ such that $\r\x+\sum^t_{i=1}\b_{i}\x_i \in L$. Since $\{\x_1,\x_2,\cdots,\x_t\}$ is
algebraically independent over $L$, we know that $\r\neq 0$. Therefore \begin{equation}\label{G action on x} \x-\sum^t_{i=1}\a_i\x_i\in L, \quad \text{where}\ \a_i:=\frac{-\b_i}{\r}\end{equation}
and thus there is an $a_\x\in F$ such that $\x=\sum^t_{i=1}\a_i\x_i+a_\x$. Note that every element of the form $\sum^t_{i=1}\a_i\x_i+a$, where $(\a_1,\cdots,\a_t)\in C^t$ and $a\in F$, is clearly an antiderivative of $F$.
\epf

%%%%%%%%%%%%%%%%%%%%%%%%%%%%%%%%%%%%%%%%%%%%%%%%%%%%%%%%%%%%%%%%%%%%%%%%%%%%%%%%%%%%%%%%%%%%%%%%%%%%%%%%%%%%%%%%
%%%%%%%%%%%%%%%%%%%%%%%%%%%%%%%%%%%%%%%%%%%%%%%%%%%%%%%%%%%%%%%%%%%%%%%%%%%%%%%%%%%%%%%%%%%%%%%%%%%%%%%%%%%%%%%
\bt \label{Structure AE} Let $E=F(\x_1,\x_2,\cdots,\x_t)$ be an antiderivative extension of $F$ and let $K$
be a differential subfield of $E$ containing $F$. Then $K$ is an antiderivative extension of $F$.\et

\bpf Let $W:=$span$_C\{\x_1,\cdots,\x_t\}$ denote the vector space generated by the elements $\x_1,\cdots,\x_t$ over the field of constants $C$ of $F$.
Let $V:=K\cap W$ and note that $V$ is a subspace of $W$. Let  $S_1\subset W$  be a $C-$ basis for  $V$. We claim that $K=F(S_1)$. Choose a set $S_2\subset W$ so that $S_1\cup S_2$ is a $C-$basis for $W$.  Clearly, $S_1\cup S_2$ is a finite set consisting of antiderivatives of $F$, the field $F(S_1)$ is  a differential field and $K\supseteq F(S_1)\supset V.$   Also note that $F(S_1\cup S_2)=F(W)=E$. If elements of $S_2$ are algebraically dependent over $K$ then by theorem \ref{ostrowski}, $K$ contains a non zero $C-$linear combination of elements of $S_2$. But then, such  a linear combination should be in $V$, a contradiction to the fact that $S_1\cup S_2$ is linearly independent over $C$.  Thus $S_2$ is algebraically independent over $K$.  Therefore,  tr.d.$(E|K)=$ tr.d.$(E|F(S_1))$ and since $K\supseteq F(S_1)$, we see that $K$ is algebraic over $F(S_1)$. Now by theorem \ref{no alg extn thm}, we obtain $K=F(S_1)$. Hence our claim. Now since $S_1\subset W$, we see that $S_1$ consists of antiderivatives of $F$ and thus $K$ is an antiderivative extension of $F$.
\epf
%%%%%%%%%%%%%%%%%%%%%%%%%%%%%%%%%%%%%%%%%%%%%%%%%%%%%%%%

%%%%%%%%%%%%%%%%%%%%%%%%%%%%%%%%%%%%%%%%%%%%%%%%%%%%%%%%
\subsection{Differential Automorphisms of Antiderivative Extensions} Let $E= F(\x_1,$ $\cdots,\x_t)$  be an antiderivative extension of $F$.  By definition, $E$ is a no new constant extension of $F$. In light of theorem \ref{alg dep of int}, we may assume $\x_1,\x_2,\cdots,\x_t$ are algebraically independent over $F$. Let $R:= F[\x_1,$ $\cdots,\x_t]\subset E$ and note that $R$ is a differential ring. Let
$\sigma\in G:=G( E| F)$. Then since $\x'_i\in F$, we have
$\sigma(\x_i)'=\sigma(\x'_i)=\x'_i.$
That is,
$\big(\sigma(\x_i)-\x_i\big)'=0$. Since $ E$ is a no new constants extension
of $ F$, there is an element $\alpha_{i\sigma}\in C$ such that
$\sigma(\x_i)-\x_i=\alpha_{i\sigma}$ and therefore,
$\sigma(\x_i)=\x_i+\alpha_{i\sigma}$. Also note that $\s(\phi(\x_i))= \x_i+\alpha_{i\s}+\alpha_{i\phi}$ $=\x_i+\alpha_{i\phi}+\alpha_{i\s}$ $=\phi(\s(\x_i))$.  Since any automorphism
of $ E$ fixing $ F$ is completely determined by its action on
$\x_1,\x_2,\cdots,\x_t,$ we see that the group $G$ is commutative and that there is an injective group homomorphism from $G$ to $(C^t,+)$ given by $\sigma\hookrightarrow(\alpha_{1\sigma},\cdots,\alpha_{t\sigma})$. To prove surjectivity, let $\vec{\alpha}:=(\alpha_1,\alpha_2,
\cdots,\alpha_n)\in C^t$. Define a ring ($F-$algebra) homomorphism $\s_{\vec{\alpha}}:R\to R$  by setting $\s_{\vec{\alpha}}(\x_i)=\x_i+\alpha_i$ and $\s_{\vec{\alpha}}(f)=f$ for all $f\in F$. The ring homomorphism obtained by mapping $\x_i\mapsto \x_i-\alpha_i$ and fixing elements of $F$ is the inverse of $\s_{\vec{\alpha}}$ and therefore $\s_{\vec{\alpha}}$ is a ring automorphism. Since $\s_{\vec{\alpha}}(\x_i)'=\s_{\vec{\alpha}}(\x'_i)$, we see that $\s_{\vec{\alpha}}$ is a differential ring automorphism. Now we extend
$\s_{\vec{\alpha}}$ to the field of fractions $E$ of $R$ to obtain a differential field automorphism.
Thus $G$ is isomorphic to the commutative group $(C^t,+)$. We refer the reader to \cite{Magid 1} and \cite{Marius-Singer} for a thorough treatment of differential fields and Picard-Vessiot theory.
%%%%%%%%%%%%%%%%%%%%%%%%%%%%%%%%%%%%%%%%%%%%%%%%%%%%%%%%%%%%%%%%%%%%%%%%%%%%
%%%%%%%%%%%%%%%%%%%%%%%%%%%%%%%%%%%%%%%%%%%%%%%%%%%%%%%%%%%%%%%%%%%%%%%%%%%
\bprop\label{fixed field}
Let $E=F(\x_1,\x_2,\cdots,\x_t)$ be an antiderivative extension of $F$. Then the  fixed field  $E^{G(E|F)}:=\{y\in E| \s(y)=y, \ \text{for all}\ \s\in G(E|F)\}$  equals $F$.

\eprop

\bpf

Let $u\in E-F$ and consider $F\langle u\rangle$, the differential field generated by $F$ and $u$.
Then by theorem \ref{Structure AE}, $F\langle u\rangle$ contains an element of the form $\sum^t_{i=1} \alpha_i\x_i$,
where at least one of the $ \alpha_i$ is non zero, say $ \alpha_1\neq 0$. Let $\vec{e}_1:=(1,0,\cdots,0)\in C^t$.
The differential automorphism $\s_{\vec{e}_1}$ induced by $\vec{e}_1$ fixes all $\x_i$ when $i\geq 2$
and maps $\x_1$ to $\x_1+1$. Therefore $\s_{\vec{e}_1}(\sum^t_{i=1} \alpha_i\x_i)\neq \sum^t_{i=1} \alpha_i\x_i$. And
since $\sum^t_{i=1} \alpha_i\x_i\in F\langle u\rangle$, we obtain $\s_{\vec{e}_1}(u)\neq u$. Thus  $E^{G(E|F)}=F$.\epf
%%%%%%%%%%%%%%%%%%%%%%%%%%%%%%%%%%%%%%%%%%%%%%%%%%%%%%%%%%%%%%%%%%%%%%%%%%%%%%%
%%%%%%%%%%%%%%%%%%%%%%%%%%%%%%%%%%%%%%%%%%%%%%%%%%%%%%%%%%%%%%%%%%%%%%%%%

%%%%%%%%%%%%%%%%%%%%%%%%%%%%%%%%%%%%%%%%%%%%%%%%%%%%%%%%%%%%%%%%%%%%%%%%%%%%%%%%%%%%%%%%%%%
%%%%%%%%%%%%%%%%%%%%%%%%%%%%%%%%%%%%%%%%%%%%%%%%%%%%%%%%%%%%%%%%%%%%%%%%%%%%%%%%%%%%%%%%%%%
\section{preparation for a structure theorem}
%%%%%%%%%%%%%%%%%%%%%%%%%%%%%%%%%%%%%%%%%%%%%%%%%%%%%%%%%%%%%%%%%%%%%%%%%%%%%%%%%%%%%%%%%%%
%%%%%%%%%%%%%%%%%%%%%%%%%%%%%%%%%%%%%%%%%%%%%%%%%%%%%%%%%%%%%%%%%%%%%%%%%%%%%%%%%%%%%%%%%%%
%%%%%%%%%%%%%%%%%%%%%%%%%%%%%%%%%%%%%%%%%%%%%%%%%%%%%%%%%%%%%%%%%%%%%%%%%%%%%%%%%%%%%%%%%%%%%%%%%%%%%%%%%%%%%%%%%%%%%%%%%%%%%%%%%%%%%%%%%%%%%%%%%%%%%%%
%%%%%%%%%%%%%%%%%%%%%%%%%%%%%%%%%%%%%%%%%%%%%%%%%%%%%%%%%%%%%%%%%%%%%%%%%%%%%%%%%%%%%%%%%%%%%%%%%%%%%%%%%%%%%%%%%%%%%%%%%%%%%%%%%%%%%%%%%%%%%%%%%%%%%%%
Hereafter, we will assume that the field of constants $C$ of $F$ is an algebraically closed field.
\subsection{Normal Tower}\label{normal tower}
%%%%%%%%%%%%%%%%%%%%%%%%%%%%%%%%%%%%%%%%%%%%%%%%%%%%%%%%%%%%%%%%%%%%%%%%%%%%%%%%%%%%%%%%%%
%%%%%%%%%%%%%%%%%%%%%%%%%%%%%%%%%%%%%%%%%%%%%%%%%%%%%%%%%%%%%%%%%%%%%%%%%%%%%%%%%%%%%%%%%%

Let $N$ be a no new constants extension of $F$. We say that $K$ is the
\textsl{antiderivative closure}
\textsl{of} $F$ \textsl{in} $N$ if $K$ is generated over $F$ by
all antiderivatives of $F$ that are in $N$. Let
$E=F(\x_1,\x_2,\cdots,\x_t)$ be an iterated antiderivative extension of $F$ and for every integer $i\geq 1$, let
$E_i$ denote the antiderivative closure of $E_{i-1}$ in $E$, where
$E_0:=F$. Since $\x_i\in E_i$, we see that $E_{t}=E$. Choose the smallest
integer $m$  such that $E_m=E_{m+1}$. Clearly such an $m$ exists,
$E_i\supset E_{i-1}$ for all $1\leq i\leq m$ and $E=E_m$. We will call the tower
\begin{equation}E=E_m\supset E_{m-1}\supset\cdots\supset
E_1\supset E_0=F\end{equation} \textsl{the normal
tower} of $E$.

We will now show that the normal tower of $E$ is kept invariant under the action of $G:=G(E|F)$. We use the notation $GK$ to denote the differential field $\{\s(y)|\s\in G\ \text{and} \ y\in K\ \}$. Since $G$ fixes $K$ and $K\supseteq F$, $G$ fixes
$E_0:=F$ and thus $GE_0\subseteq E_0$. Assume that $GE_{i-1}\subseteq E_{i-1}$ for some $i$ and
let $\y\in E_i$ be an antiderivative of $E_{i-1}$. Observe that $\s(\y)'=\s(\y')$ and since $\y'\in E_{i-1}$,
by our assumption, $\s(\y')\in E_{i-1}$. Thus, for each $\s\in G$,  $\s(\y)$ is an antiderivative of $E_{i-1}$
and therefore $\s(\y)\in E_i$. Since $E_i$ is generated as a field by antiderivatives of $E_{i-1}$,
$GE_i\subseteq E_i$. Hence by induction, $GE_i\subseteq E_i$  for all $i$.

Let $N$ be a no new constants extension of $F$. Let $\y_1,\y_2,\cdots,\y_n\in N$ be iterated antiderivatives (respectively, antiderivatives) of $F$ and let $H\subseteq G(N|F)$ be a set consisting of commuting differential automorphisms. We say the $\y_1,\y_2,\cdots,\y_n\in N$ are $H-$invariant iterated antiderivatives (respectively, $H-$invariant antiderivatives) of $F$ if $\y_1,\y_2,\cdots,\y_n$ are algebraically independent iterated antiderivatives (respectively, antiderivatives) of $F$ and for each $i$, $HF_i\subseteq F_i$, where $F_i:=F(\y_1,\y_2,\cdots,\y_{i-1})$ and $F_0:=F$.

%%%%%%%%%%%%%%%%%%%%%%%%%%%%%%%%%%%%%%%%%%%%%%%%%%%%
%%%%%%%%%%%%%%%%%%%%%%%%%%%%%%%%%%%%%%%%%%%%%%%%%%%%
\begin{example}\label{log example}
Consider the fields $ L:=\C(z,\log z,\log(\log z))$ and $\mathfrak{ L}:=\C(z, S, \mathfrak{S})$, where $S:=\{\log(z+\a)|\a\in\C\}$ and $\mathfrak{S}:=\{\log(\b+\log(z+\a))| \a,\b\in C\}$. It can be shown that $\mathfrak{ L}$ is a no new constants extension of $\C$ with respect to the usual derivation $d/dz$ and that the set $\{z\}\cup S\cup\mathfrak{S}$ consists of elements algebraically independent over $\C$, see \cite{Ravi.S}.

For convenience, we will use $'$ to denote $d/dz$. Let $K\neq \C$ be a differential subfield of $ L$.  If tr.d.$(K|\C)=3$ then since tr.d.$( L|\C)=3$, by theorem \ref{no alg extn thm} we have $K= L$. Assume  tr.d.$(K|\C)=2$. We claim that $K=\C(z,\log z)$. It is enough to show that $z,\log z \in K$. Suppose that $z\notin K$.  Then tr.d.$(K(z)|\C)=3$ and thus $K(z)= L$. Now let $\s_1\in G(K(z)|K)$ be a differential automorphism that sends $z$ to $z +1$. Since $\log z\in K(z)$ and $(\log z)'=\frac{1}{z}$, we see that  $(\s^n_1(\log z))'=\frac{1}{z+n}$, for any integer $n\geq 1$. Since $\mathfrak{ L}$ is a no new constants extension of $\C$ and $(\log(z+n))'=\frac{1}{z+n}$,  we obtain that $\log(z+n)=\s^n_1(\log z)+c_n\in  L$ for some constants $c_n\in \C$. Since the set $S$ is algebraically independent over $\C$, we obtain a contradiction to the fact that $ L$ has a finite transcendence degree over $\C$. Thus $z\in K$.

Note that if $\log z\notin K$ then $K(\log z)= L$ and there is a $\s_1\in G(K(\log z)|K)$ that sends $\log z$ to $1+\log z$. Then $\log(n+\log z)=\s^n_1(\log(\log z))+c_n\in  L$ for some $c_n\in \C$, which again contradicts the fact that $ L$ has a finite transcendence degree over $\C$. Hence the claim. Similarly, one proves that if tr.d.$K|\C=1$ then $K=\C(z)$. Thus we have shown that the differential subfields of $ L$ that contains $\C$ are $ L, \C(z,\log z), \C(z)$ and $\C$.  In deed,  the normal tower of $ L$ is
$$ L\supset \C(z,\log z)\supset \C(z)\supset \C.$$\end{example}

\begin{remark}
From the above discussion, we see that $ L$ cannot be a subfield of (or not imbeddable in) any Picard-Vessiot extension of $\C(z)$ since a Picard Vessiot extension has a finite transcendence degree over its ground field. One can list all the finitely differentially generated subfields of $\mathfrak{ L}$, see \cite{Ravi.S}. Rest of this section discusses the action of differential automorphisms on iterated antiderivatives.\end{remark}

%%%%%%%%%%%%%%%%%%%%%%%%%%%%%%%%%%%%%%%%%%%%%%%%%%%%
%%%%%%%%%%%%%%%%%%%%%%%%%%%%%%%%%%%%%%%%%%%%%%%%%%%%
\bl\label{comm end} Let $N$ be a no new constants extension of $F$ and let $E$ and $L$ be differential fields such that
$N\supseteq E\supset L\supseteq F$. Let $H$ be a commutative subset of $G(N|F)$ such that $HE\subseteq E$ and $HL\subseteq L$.  If  $E$ is an
antiderivative extension of $L$ then there are $H-$invariant antiderivatives $\y_1,\y_2,\cdots,\y_t$  of $L$ such that $E=L(\y_1,\y_2,\cdots,\y_t)$.  Moreover, for each $i$ and for each $\s\in H$, $$\s(\y_i)= \de_{i\s}\y_i+\sum^{i-1}_{j=1}\r_{ij\s}\y_j+ a_{i\s},$$
for some $ \de_{i\s},\r_{ij\s}\in C$ and $a_{i\s}\in L.$ In particular, $\s(\y_i)- \de_{i\s}\y_i\in L_{i-1}$.
\el

\bpf Suppose that $E=L(\x_1,\x_2,\cdots,\x_t)$ is an
antiderivative extension of $L$.
Since $H$ keeps $L$ and $E$ invariant, for each $\s\in G$,  $\s(\x_i)\in E$ is an antiderivative of $L$. For each $i$, we apply proposition \ref{alg dep of int} and obtain constants  $\a_{ij\s}\in C$, not all zero,
such that \begin{equation}\label{G action on x} \s(\x_i)-\sum^t_{j=1}\a_{ij\s}\x_i\in L.\end{equation}
We view the quotient space $E/L$ as a $C-$vector space (infinite dimensional)
and denote its element by $\overline{y}$,
where $y\in E$. There is natural action of $H$ on $E/L$, namely, $\s\cdot\overline{y}=\overline{\sigma(y)}$. This action is well defined since $H$ keeps $L$ and $E$ invariant.  From equation \ref{G action on x}
we see that
\begin{equation}\label{G action on y bar}\s\cdot\bx_i=\sum^t_{j=1}\a_{ij\s}\bx_i\end{equation}
for every $\s\in H$. Thus, the finite dimensional subspace $W:=$ span$_C$ $\{\bx_1,\cdots,\bx_t\}$ of $E/L$ is kept invariant
under the action of $H$. The above equation induces a group homomorphism $\Phi:H\to End(W)$ and since $H$ is commutative, $\Phi(H)$ is commutative as well. It is a well known
fact that any commuting set of endomorphisms of a vector space over an algebraically closed field\footnote{Here we use the assumption that the field of constants $C$ of $F$ is algebraically closed.} can be triangularized (see \cite{Humphreys}, page 100). That is, there is a
basis $\{\by_1,\by_2,\cdots,\by_t\}$ of $W$ and  there are constants $\r_{ij\s}\in C$ such that
\begin{equation}\label{G action on the y}
\s\cdot\by_i=\overline{\sigma(\y_i)}=\sum^{i}_{j=1}\r_{ij\s}\by_j.\end{equation}
For each $i$, we have
$\by_i=\sum^m_{j=1}\b_{ij}\bx_j$ and
therefore there are elements $r_i\in L$ such that
$\y_i=\sum^m_{j=1}\b_{ij}\x_j+ r_i.$
Thus, from proposition \ref{alg dep of int}, each $\y_i$ is an antiderivative of $L$. The linear
independence of $\{\by_i|1\leq i\leq t\}$ over $C$ and theorem \ref{ostrowski} together will guarantee the algebraic
independence of $\{\y_i|1\leq i\leq t\}$ over $L$. Since $L(\y_1,\cdots,\y_t)\subseteq E$ and tr.d.$(E|L)=$ tr.d.$(L(\y_1,\cdots,\y_t)|L)$, we may apply
theorem \ref{no alg extn thm} and obtain $E=K$. For each $i$, we set
$L_i:=L(\y_1,\cdots,\y_i)$ and observe from equation \ref{G action on the y} that $HL_i\subseteq L_i$. From equation \ref{G action on the y}, we see that $\s\y_i-\r_{ii\s}\y_i-\sum^{i-1}_{j=1}\r_{ij\s}\y_j= a_{i\s}$ for some $a_{i\s}\in L$. Thus $\s\y_i= \de_{i\s}\y_i+\sum^{i-1}_{j=1}\r_{ij\s}\y_j+ a_{i\s}$, where $ \de_{i\s}:=\r_{ii\s}$. Clearly, $\s\y_i-\de_{i\s}\y_i\in L_{i-1}$. \epf

%%%%%%%%%%%%%%%%%%%%%%%%%%%%%%%%%%%%%%%%%%%%%%%%%%%%%%%%%%%%%%%%%%%%%%%%%%%%%%%%%%%%%%%%%%%%%%%%%%%%%%%%%%%%%%%%%%%%%%%%%%%%%%%%%%%%%%%%%%%%%%%%%%%%%%%
%%%%%%%%%%%%%%%%%%%%%%%%%%%%%%%%%%%%%%%%%%%%%%%%%%%%%%%%%%%%%%%%%%%%%%%%%%%%%%%%%%%%%%%%%%%%%%%%%%%%%%%%%%%%%%%%%%%%%%%%%%%%%%%%%%%%%%%%%%%%%%%%%%%%%%%
%%%%%%%%%%%%%%%%%%%%%%%%%%%%%%%%%%%%%%%%%%%%%%%%%%%%%%%%%%%%%%%%%%%%%%%%%%%%%%%%%%%%%%%%%%%%%%%%%%%%%%%%%%%%%%%%%%%%%%%%%%%%%%%%%%%%%%%%%%%%%%%%%%%%%%%

\begin{corollary}\label{calculation 1} Let $E$ be an iterated antiderivative extension of $F$ and let $H$ be a commutative subset of $G(E|F)$. Then  there are  $H-$invariant iterated antiderivatives $\y_1,\y_2,\cdots,\y_t$  of $F$ such that $E=F(\y_1,\y_2,\cdots,\y_t)$.
Moreover, for each $i$ and each $\s\in G$, $$\s(\y_{i})= \de_{i\s}\y_{i}+r_{i\s},$$
for some $ \de_{i\s}\in C$ and $r_{i\s}\in L_{i-1}.$

\end{corollary}

\bpf
Let $E=E_m\supset E_{m-1}\supset\cdots\supset
E_1\supset E_0=F$ be the normal tower of $F$. Note that $E_j$ is an antiderivative extension of $E_{j-1}$ and from section \ref{normal tower} we know that $HE_j\subseteq E_j$ for each $j$. Thus applying lemma \ref{comm end} with $M:=E_j$ and $L:=E_{j-1}$, we obtain elements $\y_{ji}$ and $H-$invariant differential fields $L_{ji}$ for $i=1,2,\cdots,t_j$. Now we rename $\y_{11}, \cdots,\y_{1t_1},\cdots, \y_{m1},\cdots,\y_{mt_m}$ as $\y_1,\cdots,\y_t$ and $L_{11}, \cdots,L_{1t_1},\cdots, L_{m1}$ $,\cdots,L_{mt_m}$ as $ L_1,\cdots,L_t$, where $t:=\sum^m_{i=1}t_i$. One can easily check that $L_i$ and $\y_i$ satisfy the desired properties.
\epf
%%%%%%%%%%%%%%%%%%%%%%%%%%%%%%%%%%%%%%%%%%%%%%%%%%%%%%%%%%%%%%%%%%%%%%%%%%%%%%%%%%%%%%%%%%%%%%%%%%%%%%%%%%%%%%%%%%%%%%%%%%%%%%%%%%%%%%%%%%%%%%%%%%%%%%%
%%%%%%%%%%%%%%%%%%%%%%%%%%%%%%%%%%%%%%%%%%%%%%%%%%%%%%%%%%%%%%%%%%%%%%%%%%%%%%%%%%%%%%%%%%%%%%%%%%%%%%%%%%%%%%%%%%%%%%%%%%%%%%%%%%%%%%%%%%%%%%%%%%%%%%%
%%%%%%%%%%%%%%%%%%%%%%%%%%%%%%%%%%%%%%%%%%%%%%%%%%%%%%%%%%%%%%%%%%%%%%%%%%%%%%%%%%%%%%%%%%%%%%%%%%%%%%%%%%%%%%%%%%%%%%%%%%%%%%%%%%%%%%%%%%%%%%%%%%%%%%%
We need the following technical (rather computational) lemma to prove theorem \ref{structure IAE}.

\bl \label{calculation 2} Let $E$ be an iterated antiderivative extension of $F$. Suppose that $K\supseteq F$ be differential subfield of $E$ such that $E$ is an antiderivative extension of $K$ and let $G:= G(E|K)$. Then, there are  $G-$invariant iterated antiderivatives $\y_1,\y_2,$ $\cdots,\y_t$ of $F$ such that $E=F(\y_1,\cdots,\y_t)$. Let $L^*:=F(\y_1,\cdots,\y_{t-1})$. Then, either  $K\subseteq L^*$ or there is an element $a\in L^*$  such that $\y_t+a\in K$. Moreover,  $\y_t+a\notin F\langle \y'_t+a'\rangle$ and thus $F\langle \y'_t+a'\rangle$ is a proper differential subfield of $K$.\el

\bpf Since $G$ is a commutative group, from corollary \ref{calculation 1}, it follows that there are  $G-$invariant iterated antiderivatives $\y_1,\y_2,\cdots,\y_t$ of $F$ such that $E=F(\y_1,\cdots,\y_t)$. Assume that $K\nsubseteq L^*:=F(\y_1,\cdots,\y_{t-1})$ and let $u\in K\cap (E-L^*) $. Since $E=L^*(\y_t)$,
we may write $u=P/Q,$ where
$P,Q\in L^*[\y_t]$, $ P, Q$ relatively prime, and $Q$ is
monic. From corollary \ref{calculation 1}, we have
\begin{equation}\label{cruc sigma eqn of
y}\s(\y_t)= \de_{\s}\y_t+r_{\s}\end{equation}
for every $\sigma\in
G$, where $ \de_{\s}\in C$ and $r_{\s}\in L^*$. Thus $G$ consists of differential
automorphisms of the ring $L^*[\y_t]$. Since $u\in K$, we have $\s(u)=u$ for all
$\s\in G$. Thus
$\s(P)Q=\s(Q)P.$ Since $P$ and $Q$ are relatively prime, $P$
divides $\s(P)$ and $Q$ divides $\s(Q)$. But from equation \ref{cruc sigma eqn of y}, we see that deg $\s(P)=$
deg $P$ and
 deg $\s(Q)=$ deg $Q$ and thus $\s(P)=f_\s P$ and $\s(Q)=g_\s Q$ for some
$f_\s, g_\s\in L^*$. Since $\s(P/Q)=P/Q$, we must have $f_\s=g_\s$.
Now writing $Q=\sum^l_{i=0}b_i\y_t^i$ with $b_i\in L^*$ (note that $b_l=1$), we observe that
$$\sum^l_{i=0}\s(b_i)( \de_{\s}\y_t+r_{\s})^i=f_\s(\sum^l_{i=0}b_i\y_t^i).$$
Thus comparing the coefficients of $\y_t^l$, we obtain $ \de^l_{\s}=f_\s$.
Hence, for all $\s\in G$, $\s(P)= \de^l_{\s}P$ and
$\s(Q)= \de^l_{\s}Q$, where $ \de^l_{\s}\in C$. Then
$P'/P, Q'/Q\in E^G$--the fixed field of the group
$G$. From proposition \ref{fixed field}, we know that $E^G=K$ and thus $P'/P,
Q'/Q\in K$, where  $P, Q\in E$. Now from theorem \ref{no dlog} we obtain that $P,Q\in K$. Hence $G$ fixes both $P$ and $Q$.

Since $u\notin L^*$, we have $P$ or $Q$ does not belong to $L^*$. Without loss of generality, assume $P\notin L^*$. Then there is an $n\geq 1$ and $a_i\in L^*$ such that $P=\sum^n_{i=0}a_i\y_t^i$. Now, for any $\s\in G$, we have
$\s(P)=P$ and therefore
\begin{align*}&\s(a_n)( \de_{\s}\y_t+r_{\s})^n+\s(a_{n-1})( \de_{\s}\y_t+r_{\s})^{n-1}+\cdots+
\s(a_0)\\&=a_n\y_t^n+a_{n-1}\y_t^{n-1}+\cdots+a_0\end{align*}
Comparing the coefficients of $\y_t^n$, and respectively of $\y_t^{n-1}$,  we obtain
\begin{align}\label{the sigma eqn I} &\s(a_n)= \de^{-n}_{\s}a_n\quad \text{and} \\
&\label{the sigma eqn II}n \de^{n-1}_{\s}\s(a_n)r_{\s}+ \de^{n-1}_{\s}\s(a_{n-1})=a_{n-1},\end{align}

for every $\s\in G$. Since $ \de_{\s}\in C$, from equation \ref{the sigma eqn I}, we have $a'_n/a_n\in E^G=K$ and therefore applying theorem \ref{no dlog}, we obtain $a_n\in K$. In particular
$ \de^{n}_{\s}=1$. Now from equation \ref{the sigma eqn II}, we obtain
\begin{align}
&\s(a_{n-1})= \de_{\s}(a_{n-1})-n a_n r_{\s}\notag\quad \text{and thus}\\
&\s\left(a_{n-1}/n a_n\right)= \de_{\s}\left(a_{n-1}/n a_n\right)-r_{\s}\label{cruc sigma eqn alt}.
\end{align}

We add equations \ref{cruc sigma eqn alt} and \ref{cruc sigma eqn of y} to get
\begin{equation}\s\left(\y_t+\frac{a_{n-1}}{n a_n}\right)= \de_{\s}\left(\y_t+\frac{a_{n-1}}{n a_n}\right)\quad\text{for all}\ \s\in G.\end{equation}

Let $a:=a_{n-1}/n a_n$ and observe that
$(\y_t+a)'/(\y_t+a)\in E^G=K$. Again by theorem \ref{no dlog} we should then have $\y_t+a\in K$.  Note that $\y_t'+a'\in L^* $ and thus $F\langle \y_t'+a'\rangle\subseteq L^*$. And since $\y_t\notin L^*$ and $a\in L^*$ we know that $\y_t+a\notin F\langle \y_t'+a'\rangle$.  Thus $\y_t+a\in K-F\langle \y_t'+a'\rangle$ is an antiderivative of $F\langle \y_t'+a'\rangle$. Thus $\y_t+a$ is transcendental over $F\langle \y_t'+a'\rangle$ and therefore tr.d.$(K|F\langle \y_t'+a'\rangle) \geq 1$. Hence $F\langle \y_t'+a'\rangle$ is a proper differential subfield of $K$.\epf

%%%%%%%%%%%%%%%%%%%%%%%%%%%%%%%%%%%%%%%%%%%%%%%%%%%%%%%%%%%%%%%%%%%%%%%%%%%%%%%%%%%%%%%%%%%%%%%%%%%%%%%%%%%%%%%%%%%%%%%%%%%%%%%%%%%%%%%%%%%%%%%%%%%%%%%
%%%%%%%%%%%%%%%%%%%%%%%%%%%%%%%%%%%%%%%%%%%%%%%%%%%%%%%%%%%%%%%%%%%%%%%%%%%%%%%%%%%%%%%%%%%%%%%%%%%%%%%%%%%%%%%%%%%%%%%%%%%%%%%%%%%%%%%%%%%%%%%%%%%%%%%
%%%%%%%%%%%%%%%%%%%%%%%%%%%%%%%%%%%%%%%%%%%%%%%%%%%%%%%%%%%%%%%%%%%%%%%%%%%%%%%%%%%%%%%%%%%%%%%%%%%%%%%%%%%%%%%%%%%%%%%%%%%%%%%%%%%%%%%%%%%%%%%%%%%%%%%

\section{structure theorem}

%%%%%%%%%%%%%%%%%%%%%%%%%%%%%%%%%%%%%%%%%%%%%%%%%%%%%%%%%%%%%%%%%%%%%%%%%%%%%%%%%%%%%%%%%%%%%%%%%%%%%%%%%%%%%%%%%%%%%%%%%%%%%%%%%%%%%%%%%%%%%%%%%%%%%%%
%%%%%%%%%%%%%%%%%%%%%%%%%%%%%%%%%%%%%%%%%%%%%%%%%%%%%%%%%%%%%%%%%%%%%%%%%%%%%%%%%%%%%%%%%%%%%%%%%%%%%%%%%%%%%%%%%%%%%%%%%%%%%%%%%%%%%%%%%%%%%%%%%%%%%%%
 We recall that $M$ is a \textsl{minimal} differential field extension of $F$ if $M\supseteq F$ is a differential field extension such that if $K$ is a differential subfield of $M$ and $K\supseteq F$ then $M=K$ or $M=F$.
%%%%%%%%%%%%%%%%%%%%%%%%%%%%%%%%%%%%%%%%%%%%%%%%%%%%%%%%%%%%%%%%%%%%%%%%%%%%%%%%%%%%%%%%%%%%%%%%%%%%%%%%%%%%%%%%%%%%%%%%%%%%%%%%%%%%%%%%%%%%%%%%%%%%%%%
%%%%%%%%%%%%%%%%%%%%%%%%%%%%%%%%%%%%%%%%%%%%%%%%%%%%%%%%%%%%%%%%%%%%%%%%%%%%%%%%%%%%%%%%%%%%%%%%%%%%%%%%%%%%%%%%%%%%%%%%%%%%%%%%%%%%%%%%%%%%%%%%%%%%%%%

%%%%%%%%%%%%%%%%%%%%%%%%%%%%%%%%%%%%%%%%%%%%%%%%%%%%%%%%%%%%%%%%%%%%%%%%%%%%%%%%%%%%%%%%%%%%%%%%%%%%%%%%%%%%%%%%%%%%%%%%%%%%%%%%%%%%%%%%%%%%%%%%%%%%%%%
%%%%%%%%%%%%%%%%%%%%%%%%%%%%%%%%%%%%%%%%%%%%%%%%%%%%%%%%%%%%%%%%%%%%%%%%%%%%%%%%%%%%%%%%%%%%%%%%%%%%%%%%%%%%%%%%%%%%%%%%%%%%%%%%%%%%%%%%%%%%%%%%%%%%%%%
\bprop \label{Min DF} Let $E$ be an iterated antiderivative extension of $F$. Suppose that for any containments of differential fields $F\subseteq F^*\subset M^*\subseteq E$  such that $M^*$ is a minimal differential field extension of $F^*$, there is an antiderivative $\y\in E$ of $F^*$ such that $M^*=F^*(\y)$. Then, if $K$ is a differential subfield of $E$ such that $ K\supseteq F$ then $K$ is an iterated antiderivative extension of $F$.

\eprop

\bpf
 Let  $K$ be a differential subfield of $E$ such that $E\supseteq K\supset F$.  To avoid triviality, assume $K\neq F$. Then $K$ contains a  minimal differential field extension of $F$ and therefore by assumption,  $K-F$ contains an antiderivative $\y$ of $F^*$.  Assume that $K$ contains an iterated antiderivative extension $K^*$ of $F$ such that tr.d.$(K^*|F)=t$ for some $t\geq 1$. If $K^*\neq K$ then $K$ contains a minimal differential field extension of $K^*$ and therefore, by our assumption, $K-K^*$ contains an antiderivative $\y$ of $K^*$. Thus either $K^*=K$ or $K\supseteq K^*(\y)$, $\y'\in K^*$ and  tr.d.$(K^*(\y)|F)=t+1$. Since tr.d.$(K|F)<\infty$,  we have proved the proposition.  \epf

We note that to prove theorem \ref{structure IAE}, it is necessary
and sufficient to prove that the supposition statement of proposition \ref{Min DF} is always true for any iterated antiderivative extension of $F$.
%%%%%%%%%%%%%%%%%%%%%%%%%%%%%%%%%%%%%%%%%%%%%%%%%%%%%%%%%%%%%%%%%%%%%%%%%%%%%%%%%%%%%%%%%%%%%%%%%%%%%%%%%%%%%%%%%%%%%%%%%%%%%%%%%%%%%%%%%%%%%%%%%%%%%%%
%%%%%%%%%%%%%%%%%%%%%%%%%%%%%%%%%%%%%%%%%%%%%%%%%%%%%%%%%%%%%%%%%%%%%%%%%%%%%%%%%%%%%%%%%%%%%%%%%%%%%%%%%%%%%%%%%%%%%%%%%%%%%%%%%%%%%%%%%%%%%%%%%%%%%%%
\bt \label{minimal iae} Let $E$ be an iterated antiderivative extension of $F$ and let $K$ be a minimal differential  field extension of $F$ such that $E\supseteq K\supset F$. Then $K=F(\x)$ for some antiderivative $\x\in E$ of $F$. \et
\bpf
We will use an induction on $n:=$tr.d.$E|F$ to prove this theorem. From theorem \ref{no alg extn thm}, we know that tr.d.$(K|F)\geq 1$. In particular, $n\geq 1$.

Case $n=1$: we have tr.d.$(E|F)$= tr.d.$(K|F)$=1 and $E\supseteq K$. Applying corollary \ref{trd cor}, we obtain that $E=K$.

Let $n\geq 2$ and assume that the theorem holds for iterated antiderivative extensions of transcendence degree $\leq n-1$.   Let $E=E_m\supset
E_{m-1}\supset\cdots\supset E_1\supset E_0=F$ be the normal tower of $E$. Since $E\neq F$, from corollary \ref{trd cor}, we have tr.d.$(E_1|F)>0$ and thus  $E$ is an iterated antiderivative extension of
$E_1$ with  tr.d.$(E|E_1)\leq n-1$.  Note that if $F^*\supseteq E_1$ then tr.d.$(E|F^*)\leq$ tr.d.$(E|E_1)=n-1$. Then by induction, if $M^*$ and $F^*$ are differential fields such that   $E\supseteq M^*\supset F^*\supseteq E_1$ and that $M^*$ is a minimal differential field extension of $F^*$  then $M^*=F^*(\y)$ for some antiderivative $\y\in E$ of $E_1$. Therefore, by proposition \ref{Min DF}, we obtain that every differential subfield of $E$ that contains $E_1$ is an iterated antiderivative extension of $E_1$. Since $E\supseteq KE_1\supseteq E_1$, we obtain $KE_1$ is an iterated antiderivative extension of $E_1$. And since  $E_1$ is an antiderivative extension of $F$, we obtain that $KE_1$ is an iterated antiderivative extension of $F$ as well. If tr.d.$(KE_1|F)<$ tr.d.$(E|F)=n$ then  by induction, we have proved that $K$ is of the required form. Therefore we may assume
tr.d.$(KE_1|F)=$ tr.d.$(E|F)$, that is, $KE_1=E$. Then since  $E_1$ is an antiderivative extension of $F$ and $K\supset F$, we obtain that $E$ is an antiderivative extension of $K$ as well and thus $G(E|K)$ is a commutative group.

Now we apply lemma \ref{calculation 2} and obtain $G(E|K)-$invariant iterated antiderivatives $\y_1,\y_2,\cdots,\y_t$ of
$F$ such that $E=F(\y_1,\cdots,\y_t)$. If $ K\subseteq L^*:=F(\y_1,\cdots,\y_{t-1})$ then since tr.d.$(L^*|F)=$ tr.d.$(E|F) -1$ and $L^*$ is an iterated antiderivative extension of $F$, by induction, we are done. Otherwise, by lemma \ref{calculation 2}, there is an element $a\in L^*$ such that $\y_t+a\in K$, $\y_t+a\notin F\langle \y'_t+a'\rangle$ and that $F\langle\y'_t+a'\rangle$ is a proper differential subfield of $K$. Then, since $K$ is minimal extension of $F$,  $F\langle \y'_t+a'\rangle=F$. Thus we have $(\y_t+a)'=\y'_t+a'\in F$ and $\y_t+a\notin F$. Then $F(\y_t+a)$ is a differential field and $K\supseteq F(\y_t+a)\supset F$. Again, since $K$ is a minimal extension of $F$, we should have $K=F(\y_t+a)$ and by setting $\x:=\y_t+a$, we complete the proof.\epf

%%%%%%%%%%%%%%%%%%%%%%%%%%%%%%%%%%%%%%%%%%%%%%%%%%%%%%%%%%%%%%%%%%%%%%%%%%%%%%%%%%%%%%%%%%%%%%%%%%%%%%%%%%%%%%%%%%%%%%%%%%%%%%%%%%%%%%%%%%%%%%%%%%%%%%
%%%%%%%%%%%%%%%%%%%%%%%%%%%%%%%%%%%%%%%%%%%%%%%%%%%%%%%%%%%%%%%%%%%%%%%%%%%%%%%%%%%%%%%%%%%%%%%%%%%%%%%%%%%%%%%%%%%%%%%%%%%%%%%%%%%%%%%%%%%%%%%%%%%%%%
%%%%%%%%%%%%%%%%%%%%%%%%%%%%%%%%%%%%%%%%%%%%%%%%%%%%%%%%%%%%%%%%%%%%%%%%%%%%%%%%%%%%%%%%%%%%%%%%%%%%%%%%%%%%%%%%%%%%%%%%%%%%%%%%%%%%%%%%%%%%%%%%%%%%%%

\bt \label{structure IAE} Let $E$ be an iterated antiderivative extension of $F$ and let $K\supseteq F$ be a differential subfield of $E$.
Then $K$ is an  iterated antiderivative extension of $F$. \et
\bpf Follows from theorem \ref{minimal iae} and proposition \ref{Min DF}.
\epf
%%%%%%%%%%%%%%%%%%%%%%%%%%%%%%%%%%%%%%%%%%%%%%%%%%%%%%%%%%%%%%%%%%%%%%%%%%%%%%%%%%%%%%%%%%%%%%%%%%%%%%%%%%%%%%%%%%
\section{concluding remarks}\label{app} In this section we will see an application of theorem \ref{structure IAE}.
Throughout this section let $C$ be an algebraically closed field of characteristic zero and we view $C$ as a differential field with the trivial derivation. Consider the field of rational functions $C( z)$ and set $ z':=1$. Then it is easy to check that $C( z)$  is a  no new constant extension of $C$. Let $C( z)( z_1, z_2,$ $\cdots, z_t)$ be any iterated antiderivative extension of $C( z)$. We may also assume that $ z_1, z_2,$ $\cdots, z_t,$ are algebraically independent over $C( z)$. For any $u\in C( z, z_1, z_2,$ $\cdots, z_t)-C,$  theorem \ref{structure IAE} tells us the differential field $C\langle u\rangle=C(u,u',u'',\cdots)$ contains an antiderivative $\y\in C\langle u\rangle-C$ of $C$. Then,  $\y'=\a$ for some $\a\in C-\{0\}$  and we see that $\y'=(\a z)'$. Therefore, there is a $\b\in C$ such that  $\y=\a z+\b$, where $\a\in C-\{0\}$. Thus $ z\in C\langle u\rangle$.  Therefore, for each $u\in C( z, z_1, z_2,$ $\cdots, z_t)-C,$ there is an integer $n\geq 0$ and relatively prime polynomials $P, Q\in C[x_1,\cdots,x_{n+1}]$ such that \begin{equation}\label{writing x} z=\dfrac{P(u,u^{(1)},\cdots,u^{(n)})}{Q(u,u^{(1)},\cdots,u^{(n)})},\end{equation} where $u^{(i)}$ denotes the i-th derivative of $u$.

\begin{example}  Consider the differential field $\C( z,  \log z)$ with the usual derivation $d/dz$.  Then, for even a simple expression like $u:= \frac{\log z}{z}$, it can be tedious to write $ z$ in terms of $u$ and  its derivatives as in equation \ref{writing x}. In fact $ z=\dfrac{u''+uu'}{uu''-3(u')^2}.$ Since $ z_1=u z$, we see that $ z_1=\dfrac{uu''+u^2u'}{uu''-3(u')^2}$
and thus $\C\langle u\rangle=\C( z, \log z)$.\end{example}

%%%%%%%%%%%%%%%%%%%%%%%%%%%%%%%%%%%%%%%%%%%%%%%%%%%%%%%%%%%%%%%%%%%%%%%%%%%%%%%%%%%%%%%%%%%%%%%%%%%%%%%%%%%%%%%%%%
%%%%%%%%%%%%%%%%%%%%%%%%%%%%%%%%%%%%%%%%%%%%%%%%%%%%%%%%%%%%%%%%%%%%%%%%%%%%%%%%%%%%%%%%%%%%%%%%%%%%%%%%%%%%%%%%%%

%%%%%%%%%%%%%%%%%%%%%%%%%%%%%%%%%%%%%%%%%%%%%%%%%%%%%%%%%%%%%%%%%%%%%%%%%%%%%%%%%%%%%%%%%%%%%%%%%%%%%%%%%%%%%%%%%%%
\end{document}